\documentclass{article}
\usepackage{amsmath}
\usepackage{amsthm}
\usepackage{amsbsy}
\usepackage{amsopn}
\usepackage{amssymb}
\usepackage{amsfonts}
\usepackage[official ]{eurosym}

\newtheorem*{GBThm}{Gauss-Bonnet Theorem}
\newtheorem*{EGThm}{Theorema Egregium}
\title{Gaussian curvature in codimension $> 1$}
\date{}
\author{Daniel \'Alvarez-Gavela}
\begin{document}
\maketitle
\begin{abstract}
The Gaussian curvature $K$ is a fundamental geometric quantity discovered by Gauss in the case of surfaces embedded in $\mathbb{R}^3$. One can naturally extend the definition of the Gaussian curvature to arbitrary submanifolds of $\mathbb{R}^k$ so that the extrinsic interpretation of $K$, the Theorema Egregium and the Gauss-Bonnet Theorem still hold. We give a concise exposition of these classical facts.
\end{abstract}
\section{The generalized Gaussian curvature} \
Let $M$ be a closed $m$-dimensional oriented submanifold of $\mathbb{R}^k$, with $n=k-m$ the codimension. For each $p \in M$ and each unit vector $\nu \in \mathbb{R}^k$ normal to $M$ at $p$ we have a quadratic form on $T_pM$, namely $\text{II}^\nu_M(\cdot, \cdot) = < \text{II}_M( \cdot, \cdot), \nu >$. Throughout we will denote the first and second fundamental forms of $M$ by $\text{I}_M$ and $\text{II}_M$ respectively. Define $K_M^\nu= \text{det}(\text{II}^\nu_M)/\text{det}(\text{I}_M)$, where both determinants should be computed with respect to the same basis of $T_pM$. We think of $K_M^{\nu}$ as the Gaussian curvature of $M$ in the direction $\nu$. Let $NS_p \subset S^{k-1}$ denote the $(n-1)$-dimensional sphere of unit vectors normal to $M$ at $p$. The idea is to average $K_M^{\nu}$ over $NS_p$. Indeed, if $\omega_d$ is the volume of the $d$-dimensional sphere, we define the generalized Gaussian curvature of $M$ at the point $p$ to be
\[ K_M =\frac{1}{ \omega_{n-1}} \int_{NS_p}  K_M^\nu \, \, dV_{NS_p}. \]
If we agree that $\omega_0=2$, the reader can easily check that for even-dimensional hypersurfaces this definition coincides with the usual definition (namely $K_M^{\mathbf{n}}$ for the positively oriented unit normal $\mathbf{n}$). On the other hand, for odd-dimensional hypersurfaces we have $K_M=0$. In fact, odd-dimensional submanifolds of $\mathbb{R}^k$ always have $K_M=0$, irrespective of the codimension. Other than the geometric interpretation of $K_M$ as an average of ``directional Gaussian curvatures'' over all normal directions, the following two theorems justify the name given to $K_M$. \\ 
\begin{GBThm}
\[ \int_M K_M \, \, dV_M = \frac{ \omega_{k-1} }{ \omega_{n-1} } \chi(M). \]
\end{GBThm}
$\newline$
\begin{EGThm} $ \omega_{n-1} / \omega_{k-1} K_M$ is an intrinsic invariant of the metric.
\end{EGThm}
$\newline$
In fact we will see that $  \omega_{n-1} / \omega_{k-1} K_M dV_M = \text{Pff}(-\Omega/2\pi)$, where the right hand side is the Pfaffian of the matrix of curvature $2$-forms on $M$ (scaled by a factor of $-1/2\pi)$. The Chern-Gauss-Bonnet Theorem says that $\text{Pff}(-\Omega/2\pi)$ has total integral $\chi(M)$, i.e that it is a representative of the Euler class. Since the form $\text{Pff}(-\Omega/2\pi)$ is defined intrinsically in terms of the metric, it is the product of the volume form $dV_M$ and a function in $C^{\infty}(M)$ that must also be intrinsic. The content of the definition of $K_M$ is the observation that this function can be understood extrinsically as a natural generalization of the Gaussian curvature. \\ \\
The original Gauss-Bonnet formula was of course due to the two mathematicians who give name to the theorem: Gauss in the case of geodesic triangles and Bonnet in the case of a simply connected domain on a surface bounded by finitely many curves (it should be mentioned that this result was also proved independently by Binet). The global formulation of the theorem for surfaces was first stated by Walther von Dyck. Hopf generalized the result to hypersurfaces of Euclidean space by showing that the degree of the Gauss map is the Euler characteristic. This was further extended to arbitrary codimension by Fenchel and Allendoerfer independently, the latter using Weyl's work on tubes. Lipschitz, Killing and Kronecker played an important role in generalizing the Gaussian curvature to more general situations. Later Allendoerfer and Weil used the local isometric embedding theorem to prove the Gauss-Bonnet formula for arbitrary Riemannian manifolds, although it was Chern who gave the first purely intrinsic proof of the theorem in his colossal six-page paper, introducing the concept of transgression. Many different proofs and further generalizations have followed, the formula still being of interest in current research.  \\ \\
This historical account is a condensed version of Wu's exposition \cite{Historical development of the Gauss-Bonnet theorem}, which the reader can consult for further details and for an extensive list of references to the original publications. We will now prove the two theorems stated above. Thanks are in order to Sasha Zamorzaev for smart integration tricks and to Simon Brendle for inspiring coordinate computations.
\section{Proof of the Gauss-Bonnet Theorem}
Take an $\varepsilon$-neighborhood $N$ of $M$. For $\varepsilon$ small enough, $N$ is a tubular neighborhood of $M$. Its boundary $\partial N$ is a closed oriented hypersurface in $\mathbb{R}^k$, with Gauss map $g: \partial N \rightarrow S^{k-1}$. By diagonalizing $\text{II}_{\partial N}^g$ in an orthonormal basis it is easy to see that $K^g_{\partial N} =(-1)^{k-1}\det(dg)$, where we think of $dg$ as an endomorphism of subspaces of $\mathbb{R}^k$. Recall that $\chi(M)=\text{deg}(g)$. This can be seen as follows: take a tangent vector field $v$ on $M$ with isolated zeros and think of $g$ as a normal vector field on $\partial N$. Interpolate between $v$ and $g$ to obtain a vector field $w$ on $N$ that has the same zeros as $v$, with the same index at each zero. The degree of $g$ equals the sum of the degrees of $w$ restricted to small spheres around each of its zeros (where we now interpret the restricted $w$ as a map of spheres). Therefore deg$(g)$ is the sum of the indices of $v$ at its zeros, which by the Poincar\'e-Hopf Theorem is just $\chi(M)$. It follows that
\[  \omega_{k-1} \chi(M)  = \text{deg}(g) \int_{S^{k-1}} \, \, dV_{S^{k-1}} = \int_{\partial N} g^*(dV_{S^{k-1}} ) = (-1)^{k-1}\int_{\partial N} K^g_{\partial N} \, \, dV_{\partial N}. \]
The map $\pi : \partial N \rightarrow M$ that sends each $q \in \partial N$ to the point in $M$ closest to $q$ is a sphere bundle. The fibre over $p$ is $S_p=p+\varepsilon NS_p$. By the coarea formula 
\[ \omega_{k-1}  \chi(M) = (-1)^{k-1}\int_M \Big( \int_{S_p} \frac{K^g_{\partial N} }{ \text{NJ($\pi$)} } \, \, dV_{S_p} \Big) \, dV_M \]
where NJ($\pi$) is the normal Jacobian of $\pi$, i.e the determinant of $d\pi$ restricted to the orthogonal complement of its kernel. To establish the Gauss-Bonnet theorem we compute the integrand in appropriately chosen local coordinates. \\ \\
Fix $p \in M$ and $\nu \in NS_p$. We compute $K^g_{\partial N} / \text{NJ($\pi$)}$ at $p + \varepsilon \nu \in \partial N$. 
Take local coordinates $x=(x^i)$ on a neighborhood $U$ of  $p$ so that the vector fields $\partial_{x^i}=e_i$ restrict at $p$ to an orthonormal basis of $T_pM$ that diagonalizes $\text{II}_M^\nu$:
\[ < \text{II}_M(e_i, e_j), \nu > = \delta_{ij} \lambda_i. \]
Take an orthonormal basis $(\nu_j)$ of the normal space to $M$ at $p$ such that $\nu_{n}=\nu$ and extend it to a framing of the normal bundle of $M$ over $U$ by parallel transport, so that $D_{e_i} \nu_j$ is tangent to $M$ at $p$. We have an induced diffeomorphism
\[ F : U \times S^{n-1} \rightarrow \pi^{-1}(U), \qquad F(x,y) = x + \varepsilon y^i \nu_i(x). \]
Take local coordinates $\theta=(\theta^i)$ on $S^{n-1}$ near $(0, \ldots, 0, 1)$ so that at this point $\partial_{\theta^i}$ is the canonical $i$-th basis vector of $\mathbb{R}^n$ ($i <n$). Then $F$ yields local coordinates on $\partial N$ near $p + \varepsilon \nu$ via the local coordinates $(x, \theta)$ on $U \times S^{n-1}$. Note that
\[ \partial _{x^j}F = e_j + \varepsilon y^i D_{e_j}\nu_{i}, \qquad  \partial_{\theta^j} F = \varepsilon \nu_j .\]
At $p$ we have $<D_{e_j} \nu_n, e_k> = - <\nu_n , D_{e_j}e_k> = - \delta_{kj} \lambda_j$ and therefore $ \partial_{x^j} F = (1- \varepsilon \lambda_j)e_j.$ Hence the determinant of the first fundamental form at $p + \varepsilon \nu$ in these coordinates is
\[  \text{det}(\text{I}_{ \partial N} ) =  \varepsilon^{2(n-1)}\text{det}(1 - \varepsilon \text{II}^\nu_M )^2. \]
Now we consider the second fundamental form at $p + \varepsilon \nu$. We have
\[ D_{ \partial_{x^j}F } \partial_{x^i}F = \frac{\partial}{\partial x^j} ( e_i + \varepsilon y^k D_{e_i} \nu_k ) = D_{e_j}e_i + \varepsilon y^k D_{e_j}D_{e_i} \nu_k .\]
Note that $<\nu_n , \nu_n> =1$, hence $<D_{e_i} \nu_n , \nu_n> =0$, hence $< D_{e_j} D_{e_i} \nu_n, \nu_n > = - <D_{e_i} \nu_n, D_{e_j} \nu_n>$. Note also that $D_{e_j}\nu_n  = - \lambda_j e_j $. The unit normal is $g( p + \varepsilon \nu)=\nu$. Therefore at $p + \varepsilon \nu$ we have:
\[ \text{II}^\nu_{\partial N}(  \partial_{x^j}F ,  \partial_{x^i}F  ) =   <D_{ \partial_{x^j}F } \partial_{x^i}F, \nu_n>  =  \delta_{ij} \lambda_i (1- \varepsilon \lambda_i ).\]
Since $D_{e_i} \nu_j$ is tangent to $M$, $\text{II}^\nu_{\partial N} ( \partial_{x^j}F, \partial_{\theta^i}F) =0$. Finally, 
\[ \text{II}^\nu_{\partial N} ( \partial_{\theta^j}F, \partial_{\theta^i}F) = \varepsilon< \frac{\partial}{\partial \theta^j}  \nu_i , \nu_n> = - \varepsilon \delta_{ij} .\]
This shows that the determinant of the second fundamental form at $p + \varepsilon \nu$ in these coordinates is 
\[ \text{det}(\text{II}^\nu_{\partial N} )=  (-\varepsilon)^{n-1} \det(\text{II}^\nu_M ) \det( 1- \varepsilon \text{II}^\nu_M ) \]
and hence the classical Gaussian curvature at $p + \varepsilon \nu$ is 
\[ K^g_{\partial N} = \frac{ \text{det}(\text{II}^\nu_{\partial N} )}{ \text{det}( \text{I}_{\partial N} ) }=(-1)^{n-1} \frac{\det(\text{II}^\nu_M)}{ \varepsilon^{n-1} \det(1- \varepsilon \text{II}^\nu_M)  } . \]
To compute NJ($\pi$) let $P : U \times S^{n-1} \rightarrow U$ denote the obvious projection and observe that $\pi = P \circ F^{-1}$ on $\pi^{-1}(U)$, with ker($\pi$) corresponding to ker($P$) under $F$ (explicitly, ker($\pi)$ is tangent to $ S_p$). Moreover $dP = ( \text{Id} , 0 )$, $dF$ maps $T_pM$ onto the orthogonal complement of ker($\pi$) in $T_{p + \varepsilon \nu} \partial N$ and $dF(e_i) = (1-\varepsilon \lambda_i)e_i$ at the point under consideration. Finally, note that $\text{det}(\text{II}^{\nu}_M) = K^{\nu}_M$ at $p$ since there we have $\text{det}(\text{I}_M) =1$. It follows that 
\[ \text{NJ}(\pi) = \frac{1}{ \text{det}(1- \varepsilon \text{II}^\nu_M)} \quad \text{and} \quad \frac{K^g_{\partial N} }{ \text{NJ}(\pi)} =\frac{ (-1)^{n-1}}{ \varepsilon^{n-1}} K_M^\nu, \quad \text{whence by rescaling}  \]
\[  \int_{S_p} \frac{K_{\partial N} }{ \text{NJ}(\pi)}  \, \, dV_{S_p} = (-1)^{n-1}   \int_{NS_p} K_M^\nu \, \, dV_{NS_p} = (-1)^{n-1}  \omega_{n-1} K_M. \]
\[  \text{Therefore:} \quad \omega_{k-1} \chi(M)  = (-1)^{n+k}\omega_{n-1}  \int_M K_M \, \, dV_M. \]
If $m$ is odd, $\chi(M)=0$ and hence $\int_M K_M \, dV_M=0$ also (this case is silly since $K_M=0$ anyway). If $m$ is even, so is $n+k$ and hence $\omega_{k-1} \chi(M) = \omega_{n-1}\int_M \, \, dV_M$. This completes the proof of the Gauss-Bonnet formula.
\section{Proof of the Theorema Egregium}
To show that $\omega_{n-1}/\omega_{k-1}K_M dV_M = \text{Pff}( -\Omega/2\pi)$ we compute both differential forms at a point $p$ using the local coordinates obtained by parametrizing $M$ near $p$ as a graph over $T_pM$. In this way we can assume without loss of generality that $p=0$ and $M$ is the graph of a function $X=(X^i): \mathbb{R}^m \rightarrow \mathbb{R}^n$ such that $X(0)=0$ and $DX_0=0$. So $T_pM = \mathbb{R}^m \times 0 \subset \mathbb{R}^m \times \mathbb{R}^n =\mathbb{R}^k$ and our parametrization is $x \mapsto \big(x, X(x)$\big). We have
\[ \partial_{i} = ( \epsilon_i\,,\, \partial_{i}X ), \quad \text{where} \, \,  \epsilon_i \,\, \text{is the $i$-th canonical basis vector of } \mathbb{R}^m \, \, \text{and} \]
\[  D_{\partial_{i}}\partial_{j} = (0 \, , \, \partial^2_{i,j} X ).\]
Note that at the origin $(\partial_{i})$ is an orthonormal basis of $T_0M$, namely $(\epsilon_i)$. Let $\nu=(0, (\nu_1, \ldots, \nu_n)) \in 0 \times S^{n-1} = NS_0$. Then we have the following:
\[ \text{II}^\nu_M ( \partial_{i}, \partial_{j} ) = \nu_{\alpha} \partial^2_{i, j}X^{\alpha}, \]
\[ K_M^\nu = \sum_{ \sigma \in S_m} (-1)^{\sigma}  (\partial^2_{1, \sigma(1)}X^{\alpha_1} ) \cdots (\partial^2_{m, \sigma(m)}X^{\alpha_m} )\nu_{\alpha_1} \cdots \nu_{\alpha_m}, \]
\[ K_M = \frac{1}{\omega_{n-1}} \sum_{ \sigma \in S_m} (-1)^{\sigma}  (\partial^2_{1,\sigma(1)}X^{\alpha_1} ) \cdots (\partial^2_{m, \sigma(m)}X^{\alpha_m} ) \int_{S^{n-1}}\nu_{\alpha_1} \cdots \nu_{\alpha_m} dV_{S^{n-1}}. \]
If $m$ is odd, all the integrals are zero, so $K_M=0$. But of course $\text{Pff}(-\Omega/2\pi)$ only makes sense when $m$ is even. So from now on we assume that $m=2r$. In fact the integral over $S^{n-1}$ of the function $\nu_{\alpha_1} \cdots \nu_{\alpha_m}$ will be zero unless each $\nu_i$ appears in the monomial with an even exponent. If indeed $\nu_{\alpha_1} \cdots \nu_{\alpha_m} = \nu_1^{2a_1} \cdots \nu_n^{2a_n} $ (so that $a_1 + \cdots + a_n=r$) we can compute the integral as follows:
\[ \int_{0}^{\infty} r^{k-1} e^{-r^2} dr \int_{S^{n-1} }\nu_1^{2a_1} \cdots \nu_n^{2a_n} \,  dV_{S^{n-1} } = \int_{\mathbb{R}^n} x_1^{2a_1} \cdots x_n^{2a_n} e^{ - \sum x_i^2} \, dx_1 \cdots dx_n \]
\[ = \prod \int_\mathbb{R} x^{2a_i} e^{-x_i^2} dx_i= \prod \Gamma( a_i + 1/2 ), \quad \text{whence}   \]
\[  \int_{S^{n-1} }\nu_1^{2a_1} \cdots \nu_n^{2a_n} \,  dV_{S^{n-1} } = 2 \frac{ \prod  \Gamma( a_i + 1/2 ) }{ \Gamma( k/2) }\, \,  \, \, \text{and since} \, \, \, \omega_{n-1} = 2\frac{\pi^{n/2}}{ \Gamma(n/2)}, \, \, \text{we have}   \]
\[ K_M =  \frac{\Gamma(n/2)}{\Gamma(k/2) \pi^{n/2} }\sum_{ \sigma \in S_m, \alpha_t} (-1)^{\sigma}   \prod \Gamma( a_i + 1/2 )   (\partial^2_{1, \sigma(1)}X^{\alpha_1} ) \cdots (\partial^2_{m, \sigma( m)}X^{\alpha_m} ) \]
where we're summing over all choices of indices $1 \leq \alpha_1 , \ldots ,  \alpha_m \leq n$ such that $\nu_{\alpha_1} \cdots \nu_{\alpha_m} = \nu_1^{2a_1} \cdots \nu_n^{2a_n} $ holds for some partition $r=a_1 + \cdots + a_n$. Since
\[ \prod \Gamma( a_i + 1/2)  =  \frac{ \pi^{n/2}}{2^m} \prod \frac{(2a_i)!}{  a_i!} \quad \text{and} \quad \frac{ \omega_{n-1} }{\omega_{k-1} } =  \frac{ \Gamma(k/2)}{\pi^r  \Gamma(n/2) },  \] 
\[ \text{we deduce} \quad \frac{ \omega_{n-1} }{ \omega_{k-1} }K_M =  \frac{1}{2^m\pi^r}\sum_{ \sigma \in S_m, \alpha_t} (-1)^{\sigma}  \prod \frac{(2a_i)!}{  a_i!} (\partial^2_{1, \sigma(1)}X^{\alpha_1} ) \cdots (\partial^2_{m, \sigma( m)}X^{\alpha_m} ).  \]
We now turn to the Pffafian. The curvature $2$-forms $\Omega^l_k$ corresponding to our orthonormal basis $(\epsilon_i)$ of $T_0M$ are given by
\[\Omega^l_k = \frac{1}{2}R_{ijkl}\epsilon^i \wedge \epsilon^j, \]
where $(\epsilon^j)$ is the dual basis and $R_{ijkl}$ are the coefficients of the Riemann curvature tensor $R(X,Y,Z,W)=<[ \nabla_X, \nabla_Y ]Z - \nabla_{[X ,Y]}Z, W>$ with respect to this basis. By the Gauss equation
\[ R(X,Y,Z,W) = < \text{II}_M(X,W), \text{II}_M(Y, Z) > - < \text{II}_M(X,Z) , \text{II}_M(Y,W) > \]
\[\text{we deduce} \quad  R_{ijkl} = \sum_\gamma (\partial^2_{i,l} X^\gamma)( \partial^2_{j,k} X^\gamma) -  (\partial^2_{i,k}X^\gamma)( \partial^2_{j,l}  X^\gamma).  \]
The usual formula for the Pfaffian is
\[ \text{Pff}(-\Omega/2\pi) = \frac{(-1)^r}{2^m \pi^rr!}\sum_{ \tau \in S_m}  (-1)^{\tau} \Omega^{\tau( 2)}_{\tau( 1)} \wedge \cdots  \wedge\Omega^{\tau (m)}_{\tau(m-1)} \]
and by our choice of coordinates we have $dV_M = \epsilon^1 \wedge \cdots \wedge \epsilon^m$ at $0$, so that
\[ \text{Pff}(-\Omega/2\pi) =  \Big((-1)^r\sum_{\tau, \eta \in S_m}  \frac{(-1)^{\eta \tau}}{ 2^{m+r} \pi^r r!} R_{ \eta(1) \eta (2) \tau(1) \tau(2)} \cdots R_{ \eta(m-1) \eta(m) \tau(m-1) \tau(m) } \Big) dV_M. \]
We must compare the sum inside the brackets with our formula for $\omega_{n-1}/\omega_{k-1}K_M$. Using the above expression for $R_{ijkl}$ the sum becomes
\[\sum_{\tau, \eta \in S_m} \frac{ (-1)^r (-1)^{ \eta \tau}}{2^{m+r} \pi^r r!} \prod_{t=1}^r \Big(\sum_{\gamma_t}(\partial^2_{\eta(2t-1),\tau(2t)} X^{\gamma_{t}})( \partial^2_{ \eta(2t), \tau(2t-1)} X^{\gamma_{t}}) -  (\partial^2_{\eta(2t-1), \tau(2t-1)} X^{\gamma_{t}})( \partial^2_{\eta(2t), \tau(2t)} X^{\gamma_{t}})\Big)  \]
\[ =\sum_{\tau, \eta \in S_m, \gamma_t} \frac{(-1)^{\eta \tau}}{2^m \pi^r r!} (\partial^2_{ \eta(1), \tau(1)}X^{\gamma_1} )( \partial^2_{\eta(2), \tau(2)} X^{\gamma_1} ) \cdots ( \partial^2_{\eta(m-1), \tau(m-1)}X^{\gamma_r} )(\partial^2_{ \eta(m), \tau(m)} X^{\gamma_r}). \]
For each choice of indices $\gamma_1, \ldots , \gamma_r$, let $(\beta_1, \ldots , \beta_r)$ be the $r$-uple consisting of the same indices but reordered so that $1 \leq \beta_1 \leq \cdots \leq \beta_r \leq n$.  If $\nu_{\beta_1}^2 \cdots \nu_{\beta_r}^2 = \nu_1^{2a_1} \cdots \nu_n^{2a_n}$, the number of choices of $\gamma_1, \ldots , \gamma_r$ that correspond to $(\beta_1, \ldots , \beta_r)$ is $ r!/\prod a_i ! $. It follows that we can rewrite the sum as
\[ \sum_{ \tau, \eta \in S_m, (\beta_t)} \frac{ (-1)^{\eta \tau}} {2^m \pi^r \prod a_i!} (\partial^2_{ \eta(1), \tau(1)}X^{\beta_1} )( \partial^2_{\eta(2), \tau(2)} X^{\beta_1} ) \cdots ( \partial^2_{\eta(m-1), \tau(m-1)}X^{\beta_r} )(\partial^2_{ \eta(m), \tau(m)}  X^{\beta_r}). \]
Let's look back at our expression for $\omega_{n-1}/\omega_{k-1}K_M$. For each $r$-uple $(\beta_1, \ldots , \beta_r)$ such that $1 \leq \beta_1 \leq \cdots \leq \beta_r \leq n$, how many choices of indices $\alpha_1, \ldots , \alpha_m$ are there in the sum such that $\nu_{\alpha_1} \cdots \nu_{\alpha_m} = \nu_{\beta_1}^2 \cdots \nu_{\beta_r}^2$? Well, if the monomial is $\nu_1^{2a_1} \cdots \nu_n^{2a_n}$, there are exactly $m! / \prod(2a_i)!$ possible choices. Therefore:
\[  \frac{ \omega_{n-1} }{ \omega_{k-1} } K_M = \sum_{ \sigma, \kappa \in S_m, \alpha_ t}\frac{(-1)^{\kappa \sigma}}{2^{m}\pi^rm!} \prod \frac{ (2a_i)! }{a_i!}   (\partial^2_{\kappa(1), \sigma(1)}X^{\alpha_1} ) \cdots (\partial^2_{\kappa(m), \sigma( m)}X^{\alpha_m} ) \]
\[ = \sum_{ \sigma, \kappa \in S_m, (\beta_ t)}\frac{ (-1)^{\kappa \sigma} }{2^{m} \pi^r \prod a_i!}   (\partial^2_{\kappa(1), \sigma(1)}X^{\beta_1} )(\partial^2_{\kappa(2), \sigma(2)}X^{\beta_1} )  \cdots (\partial^2_{\kappa(m-1), \sigma(m-1)}X^{\beta_r} ) (\partial^2_{\kappa(m), \sigma( m)}X^{\beta_r} ).  \]
We have thus proved that $ \omega_{n-1} / \omega_{k-1} K_M dV_M = \text{Pff}( -\Omega / 2\pi)$.

\end{document}